\documentclass{amsart}

\usepackage{graphicx}
\newtheorem{defn}{Definition}
\newtheorem{rem}{Remark}
\newtheorem{algo}{Algorithm}

\begin{document}

\title[Oscillatory Integral over $\exp(i\pi x)x^{1/x}$]{Numerical Evaluation Of the Oscillatory Integral over $\exp(i\pi x)x^{1/x}$
between 1 and Infinity}

\author{Richard J. Mathar}
\urladdr{http://www.strw.leidenuniv.nl/~mathar}
\email{mathar@strw.leidenuniv.nl}
\address{Leiden Observatory, Leiden University, P.O. Box 9513, 2300 RA Leiden, The Netherlands}

\subjclass[2010]{Primary 65D30, 65B05; Secondary 40A25, 65T40}

\date{\today}
\keywords{Oscillatory Integral, Filon Quadrature, Incomplete Gamma Function}

\begin{abstract}
Real and imaginary part of the limit $2N\to \infty$ of the
integral $\int_1^{2N} \exp(i\pi x) \sqrt[x]{x} dx$
are evaluated to 20 digits with brute force methods after multiple
partial integration, or combining a standard Simpson
integration over the first half wave with series acceleration techniques for the
alternating series co-phased to each of its points. The integrand is of
the logarithmic kind; its branch cut
limits the performance of
integration techniques that
rely on smooth higher order derivatives.
\end{abstract}

\maketitle

\section{Scope}
\subsection{M. R. Burns' Constant}

\begin{defn}
The MRB constant is the sum of the series \cite[A037077]{EIS}
\begin{equation}
M\equiv \lim_{N\to \infty} \sum_{n=1}^{2N} (-1)^n \sqrt[n]{n}
=
\sum_{k=1}^\infty (-1)^k (k^{1/k}-1)
\approx 0.18785964
.
\label{eq.mrbdef}
\end{equation}
\end{defn}
Direct summation of the alternating series is slow and generates roughly
3 valid digits after ten thousand terms (Table \ref{tabl.dsum}). 
Euler summation \cite[(3.6.27)]{AS}\cite{Hardy} is successful in accelerating
the convergence, witnessed in Table \ref{tabl.esum}.
\begin{table}[h]
\caption{Partial sums of (\ref{eq.mrbdef}) as a function of the upper limit of summation.
}
\begin{tabular}{rl}
$\hat k$ &   $\sum_{k=1}^{\hat k} (-1)^k( k^{1/k}-1)$\\
\hline
$10$ &   0.313231759254\ldots \\
$10^2$ & 0.211329543346\ldots \\
$10^3$ & 0.191323989712\ldots \\
$10^4$ & 0.188320351076\ldots \\
$10^5$ & 0.187917210140\ldots \\
\hline
\end{tabular}
\label{tabl.dsum}
\end{table}

\begin{table}[h]
\caption{Approximations to (\ref{eq.mrbdef}) after Euler resummation
of the first $\hat k$ terms.
}
\begin{tabular}{rl}
$\hat k$ &   $E(1)$\\
\hline
  10& 0.187885886113800730351382438464680824292327407645248188116946\ldots \\
  20& 0.187859649854050194658445181421685949965109596411421113134589\ldots \\
  40& 0.187859642462068655529781630996634559217485607603915205816966\ldots \\
 100& 0.187859642462067120248517934054273314215151463271236583869269\ldots \\
 200& 0.187859642462067120248517934054273230055903094900138786171982\ldots \\
\hline
\end{tabular}
\label{tabl.esum}
\end{table}

More efficient methods lead to even  quicker convergence, as
demonstrated in Table \ref{tabl.csum}. An accuracy of 60 digits
is reached after 100 terms and will be sufficient for all purposes of this script.
\begin{table}[h]
\caption{Approximations to (\ref{eq.mrbdef}) with the first
Cohen-Villegas-Zagier algorithm using terms up to $\hat k$
\cite{CohenExpM9}.
}
\begin{tabular}{rl}
$\hat k$ &  $M$  \\
\hline
10&  0.187859642389333316567457476113016727048599369932998191180459\ldots \\
20&  0.187859642462067119674255755542940758484982176117045969528027\ldots \\
40&  0.187859642462067120248517934054273140023454509840554949525330\ldots \\
100& 0.187859642462067120248517934054273230055903094900138786171986\ldots \\
\hline
\end{tabular}
\label{tabl.csum}
\end{table}

\subsection{Oscillatory Integral}
\begin{figure}[hbt]
\includegraphics[scale=0.5,clip]{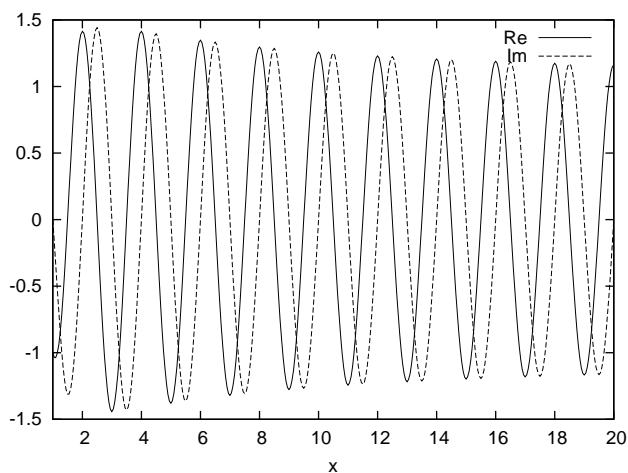}
\caption{
Real and imaginary part of $e^{i\pi x}x^{1/x}$, the integrand in (\ref{eq.ingr}).
}
\label{fig.plotI}
\end{figure}

The integrated analog of the series is a complex-valued integral of
oscillatory character,
which is difficult to evaluate by direct integration if the upper limit
becomes large, illustrated by Figure \ref{fig.plotI}.
\begin{defn} (Sequence of oscillatory integrals)
\begin{equation}
I(2N)\equiv \int_1^{2N} (-1)^x \sqrt[x]{x} dx = \int_1^{2N} e^{i\pi x}x^{1/x}dx,\quad
N\in\mathbb{Z}.
\label{eq.ingr}
\end{equation}
\end{defn}
Not convergent in the continuum limit at $N\to \infty$, the limit
of the sequence of integrals with
an integral difference in the upper limits $2N$ exists. The
objective of this work is to evaluate this limit $M_I$.
\begin{defn} (Ultraviolet limit of the sequence)
\begin{equation}
M_I \equiv \lim_{N\to\infty} I(2N)
.
\end{equation}
\end{defn}
\begin{rem}
The absolute value $|M_I|\approx 0.6876523689$ is close to $M+\frac{1}{2}$ \cite[A157852]{EIS}.
Changing the upper limit to
$2N+1$ increases $M_I$ by $2i/\pi$.
\end{rem}

To compute
$M_I$, the manuscript looks at
repeated partial integration to quench the integrand at large $x$
in preparation for standard methods of sampling along the abscissa
(Sections \ref{sec.pI} and \ref{sec.exp}),
investigates splitting the integral into an alternating series and a base interval
(Section \ref{sec.alts}),
expansion of $x^{1/x}$ into a series over $(\log x/x)^n$ (Section \ref{sec.Tlog}),
changing the path of integration in the complex plane (Section \ref{sec.cont}),
and
considers
reverse application of
the Euler-Maclaurin integral formula (Appendix \ref{app.EMcl}).

\section{Numerical Analysis} 
\subsection{Iterated Partial Integration} \label{sec.pI} 
\begin{figure}[hbt]
\includegraphics[scale=0.5,clip]{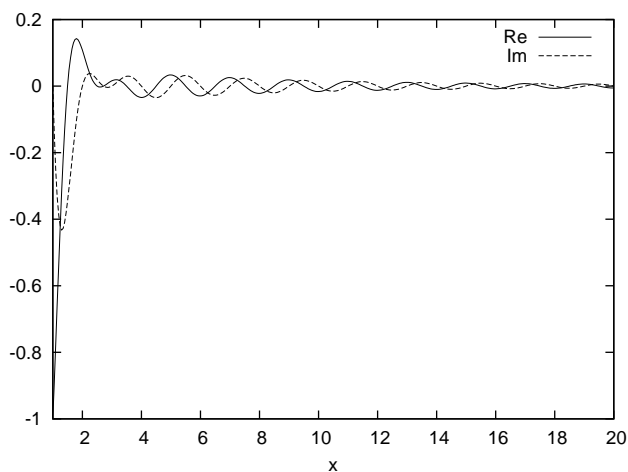}
\caption{
Real and imaginary part of the integrand in (\ref{eq.pI}).
}
\label{fig.plotIpI}
\end{figure}

A partial integration of (\ref{eq.ingr}) yields
\begin{equation}
\int_1^{2N} e^{i\pi x}x^{1/x} dx
=
-\frac{i}{\pi}e^{i\pi x} x^{1/x}\bigg|_1^{2N} +
\frac{i}{\pi} \int_1^{2N} e^{i\pi x}x^{1/x} \frac{1-\log x}{x^2}dx.
\end{equation}
The limit $N\to\infty$ can be performed in the pre-integrated term,
\begin{equation}
M_I=
-\frac{2i}{\pi}+
\frac{i}{\pi} \int_1^\infty e^{i\pi x}x^{1/x} \frac{1-\log x}{x^2}dx
\label{eq.pI},
\end{equation}
which
essentially compresses the oscillations with a factor $\propto \log(x)/x^2$,
as shown in Figure \ref{fig.plotIpI}.

A second partial integration of (\ref{eq.pI}) may follow,
\begin{multline}
\int_1^{2N} e^{i\pi x}x^{1/x} \frac{1-\log x}{x^2}dx
\\
=
-\frac{i}{\pi} e^{i\pi x} x^{1/x}\frac{1-\log x}{x^2}\Big|_1^{2N}
+\frac{i}{\pi}\int_1^{2N} e^{i\pi x}x^{1/x} \frac{1-3x+2(x-1)\log(x)+\log^2 x}{x^4}dx
\end{multline}
with the limit
\begin{equation}
M_I
=
-\frac{2i}{\pi}+
\frac{1}{\pi^2}
-\frac{1}{\pi^2}\int_1^{\infty} e^{i\pi x}x^{1/x} \frac{1-3x+2(x-1)\log(x)+\log^2 x}{x^4}dx
.
\label{eq.pI2}
\end{equation}
Repeating, a simple scheme for the $n$-th derivative of the base function $f$,
\begin{equation}
f^{(0)}(x)\equiv x^{1/x};\quad f^{(n+1)}(x)\equiv \frac{d}{dx} f^{(n)}(x),
\end{equation}
can be phrased as a set of coefficients $\alpha$,
\begin{equation}
f^{(n)}(x)\equiv x^{1/x}\sum_{r,s\ge 0} \alpha_{n,r,s}\frac{\log ^s(x)}{x^r}.
\end{equation}
Explicit computation with the chain rule establishes the recurrence
\begin{equation}
\alpha_{n+1,r,s}
=
\alpha_{n,r-2,s}
-
\alpha_{n,r-2,s-1}
+
(s+1)\alpha_{n,r-1,s+1}
-(r-1)\alpha_{n,r-1,s}
.
\end{equation}
The initial conditions are
\begin{equation}
\alpha_{n,r,s}=0\quad\text{if}\quad r<0\quad \text{or}\quad s<0;\quad \alpha_{0,0,0}=1.
\end{equation}
Equations (\ref{eq.pI}), (\ref{eq.pI2}) and the representation through $n$-fold partial integration
are summarized with
\begin{equation}
M_I = b(n) +\left(\frac{i}{\pi}\right)^n \int_1^\infty e^{i\pi x}f^{(n)}(x) dx
,
\label{eq.Mnthp}
\end{equation}
defining pre-integrated terms
\begin{equation}
b(n+1) =b(n)-\frac{i}{\pi} f^{(n-1)}(1);\quad b(0)=0;\quad b(1)=-\frac{2i}{\pi}
.
\end{equation}

The infinite interval from 1 to $\infty$ can be mapped onto the interval
from 0 to $1/2$ with the substitution (one of a family of rational maps)
\begin{equation}
y=\frac{1}{1+x};\quad x=\frac{1}{y}-1;\quad dx =-\frac{dy}{y^2};
\end{equation}
\begin{equation}
\int_1^\infty e^{i\pi x}f^{(n)}(x) dx = 
\int_0^{1/2} e^{i\pi x(y)} f^{(n)}(x(y))\frac{dy}{y^2}
.
\label{eq.map12}
\end{equation}
A combination of (\ref{eq.Mnthp}) and (\ref{eq.map12}) yields:
\begin{algo}\label{alg.pI}
Compute $M_I$ with a Simpson integration
on $s$ points in the interval $0\le y \le 1/2$:
\begin{equation}
M_I = b(n) +\left(\frac{i}{\pi}\right)^n 
\int_0^{1/2} e^{i\pi x(y)} f^{(n)}(x(y))\frac{dy}{y^2}.
\end{equation}
\end{algo}
Table \ref{tab.pI} illustrates that a choice of $n$ near 5 or 6 yields
optimum convergence (because $b(n)$ then approximate $M_I$ best),
and that integration with $s\approx 60\,000$ abscissa
points generates of the order of 13 valid digits.
Note that Romberg (Richardson) extrapolation with the standard
$15:1$ weighting of step width halving does not work
as the integrand is not in the polynomial class.
\begin{table}
\caption{Algorithm \ref{alg.pI}:
Simpson integration with $s$ points in the interval $0\le y\le 1/2$
and the $n$-th partial integration (\ref{eq.Mnthp}) inserted in (\ref{eq.map12}).
}
\begin{tabular}{rrrr}
$n$ & $s$ & $\Re(M_I)$ & $\Im(M_I)$
\\
\hline
2&2000&0.07077873792033467823&-0.68400849330911320239
\\
2&4000&0.07077826545794147869&-0.68400115445281602856
\\
\hline
3&2000&0.07077594388618778878&-0.68400039709147911789
\\
3&4000&0.07077603597753547955&-0.68400040477320427448
\\
3&8000&0.07077604728595562114&-0.68400040978725640089
\\
3&16000&0.07077604069379310873&-0.68400038906169127559
\\
3&32000&0.07077603794889411182&-0.68400038968332537849
\\
3&64000&0.07077603964902400034&-0.68400038933267435379
\\
\hline
5&2000&0.07077603934743949117&-0.68400038942675266591
\\
5&4000&0.07077603931081786932&-0.68400038943414042302
\\
5&8000&0.07077603931077043492&-0.68400038943655207372
\\
5&16000&0.07077603931149307863&-0.68400038943803823084
\\
5&32000&0.07077603931155557088&-0.68400038943794081422
\\
5&64000&0.07077603931152496616&-0.68400038943793426656
\\
\hline
6&2000&0.07077603931115051881&-0.68400038943952383435
\\
6&4000&0.07077603931156321151&-0.68400038943794556862
\\
6&8000&0.07077603931154254337&-0.68400038943792386069
\\
6&16000&0.07077603931152730519&-0.68400038943793206498
\\
6&32000&0.07077603931152869782&-0.68400038943793231267
\\
6&64000&0.07077603931152878997&-0.68400038943793211168
\\
\hline
8&2000&0.07077603932853563807&-0.68400038944555242777
\\
8&4000&0.07077603931259171964&-0.68400038943840850746
\\
\hline
9&2000&0.07077603925904942234&-0.68400038951437536286
\\
9&4000&0.07077603930824878310&-0.68400038944270983422
\\
\hline
\end{tabular}
\label{tab.pI}
\end{table}

\begin{rem}
In a variant of this mapping on a finite interval,
the transformation $x=1/u$ in the $n$-th partial integration yields
\begin{equation}
M_I=b(n)+(i/\pi)^n\int_0^1 e^{i\pi /u}f^{(n)}(1/u)\frac{du}{u^2}
.
\label{eq.1u}
\end{equation}
The precision is worse
than with Algorithm \ref{alg.pI} by approximately one digit.
I have not looked into advanced schemes for this type of oscillatory integrals
\cite{GautschiJCompAM184,HascelikJCAM223} or Sidi's generalized methods of
extrapolation.
\end{rem}

\subsection{Exponential Scaling} \label{sec.exp}
A characteristic of Algorithm \ref{alg.pI} is that the integrand
is basically reduced by another factor $1/x$ for each additional partial integration.
The variable transformation
$\log x=z$, $x=e^z$, $dx =e^zdz$,
helps to achieve exponential scaling as the integration
variable heads towards infinity, at the cost of an irregular chirp factor
in the complex exponential:
\begin{equation}
F_m=\int_m^\infty e^{i\pi x}x^{1/x}\frac{1-\log x}{x^2}dx
=\int_{\log m}^\infty (1-z)e^{i\pi \exp(z)} e^{z \exp(-z)-z} dz
.
\label{eq.F1exp}
\end{equation}
Based on
\begin{equation}
\int e^{i\pi \exp(z)+z} dz
=
-\frac{i}{\pi} e^{i\pi \exp(z)} ,
\end{equation}
and ``borrowing'' a factor $e^z$ in (\ref{eq.F1exp}) in the integrand,
\begin{equation}
F_m
=
\int_{\log m}^\infty (1-z)e^{i\pi \exp(z)+z} e^{z \exp(-z)-2z} dz,
\end{equation}
a partial integration generates
\begin{equation}
F_m
=
\frac{i(-)^m}{\pi} 
\sqrt[m]{m}
\frac{1-\log m}{m^2}
+\frac{i}{\pi}
\int_{\log m}^\infty
e^{i\pi \exp(z)}
e^{z \exp(-z)-2z}\left[2z-3+e^{-z}(1-2z+z^2)\right]dz
.
\end{equation}
Alternatively, this results applying the substitution $\log x=z$ to
(\ref{eq.Mnthp}).
After this scheme of partial integrations has been repeated $n$ times, the non-oscillating
exponential factor in the integrand is $\exp[z \exp(-z)-(n+1)z]$.
\begin{algo}\label{alg.esc}
Perform $n$ partial integrations of (\ref{eq.F1exp}), then use another
transformation $u=e^{-z}$, $z=-\log u$ to map the range $0\le z\le \infty$
to $0\le u\le 1$ in the remaining integral, and integrate
this over $u$ with a Simpson method. Eventually insert this $F_1$ in
\begin{equation}
M_I=-\frac{2i}{\pi}+\frac{i}{\pi}F_1
\end{equation}
as seen in (\ref{eq.pI}).
\end{algo}
\begin{table}
\caption{Convergence of Algorithm \ref{alg.esc}: integration
of the $n$-th partial integration of $F_1$ with $s$ equidistant
points in the interval $0\le u\le 1$.
}
\begin{tabular}{rrll}
$n$ & $s$ & $\Re(M_I)$ & $\Im(M_I)$ \\
\hline
2&4000&0.07077612979610666804&-0.68400038256040301228\\
2&8000&0.07077604264870771610&-0.68400037471801246748\\
\hline
3&4000&0.07077603954212465077&-0.68400039095784542350\\
3&8000&0.07077603950170386538&-0.68400038943798882953\\
\hline
4&4000&0.07077603927887287271&-0.68400038944795068337\\
4&8000&0.07077603931218586254&-0.68400038944147410684\\
\hline
5&4000&0.07077603931112220214&-0.68400038943707487961\\
5&8000&0.07077603931144862845&-0.68400038943795645974\\
5&16000&0.07077603931151256105&-0.68400038943794236597\\
5&32000&0.07077603931153004386&-0.68400038943793234690\\
5&64000&0.07077603931152889348&-0.68400038943793195977\\
\hline
6&4000&0.07077603931156163745&-0.68400038943790335638\\
6&8000&0.07077603931152840681&-0.68400038943792930582\\
6&16000&0.07077603931152865655&-0.68400038943793183965\\
6&32000&0.07077603931152880374&-0.68400038943793214262\\
6&64000&0.07077603931152880496&-0.68400038943793212992\\ 
6&128000&0.07077603931152880345&-0.68400038943793212926\\ 
6&256000&0.0707760393115288035480930&-0.6840003894379321291922485\\ 
6&512000&0.0707760393115288035386982&-0.6840003894379321291820339\\ 
6&1024000&0.0707760393115288035395359&-0.6840003894379321291828037\\ 
6&2048000&0.0707760393115288035395336&-0.6840003894379321291827445 \\ 
\hline
8&4000&0.07077603931135483266&-0.68400038943809843545\\
8&8000&0.07077603931151793368&-0.68400038943794252308\\
\hline
9&4000&0.07077603931088567214&-0.68400038943711147080\\
9&8000&0.07077603931148860789&-0.68400038943788083818\\
\hline
\end{tabular}
\label{tab.esc}
\end{table}
An accuracy of $10^{-21}$ can be reached sampling two million
points (Table \ref{tab.esc}) and is reported in the summary.
Comparison with Table \ref{tab.pI} shows that
one to two digits are gained relative to Algorithm \ref{alg.pI}.

\subsection{Longman's Method}\label{sec.alts}
An integral $F$ with an undulating trigonometric factor multiplied
by a monotonous
$g(x)$ may be split into an integral over the half period with
an alternating series attached to each point in that interval \cite{LynessJCAM12,LongmanCPC52,EspelidNumAlg8}.
\begin{multline}
F_m=\int_m^\infty e^{i\pi x}g(x)dx
=
\sum_{l\ge 0} \int_0^1 e^{i\pi (m+l+y)}g(m+l+y) dy
\\
=
(-)^m \int_0^1 e^{i\pi y} \sum_{l\ge 0} (-)^l g(m+l+y) dy
.
\label{eq.Fhw}
\end{multline}
The $l$-series is only alternating if the function $g(x)$ is monotonous
and does not change sign. In addition, Euler resummation assumes that
the series converges. With $M_I$, $x^{1/x}$ has a maximum at $x=e$, 
so we integrate over
$1\le x\le 3$ with any other method, setting $m=3$, and relay by one
partial integration, $g(x)=f^{(1)}(x)$, to feature a $g(x)$ that has a single
sign with decreasing $|g'(x)|$ in $[m,\infty)$.
\begin{algo} \label{alg.as}
Compute $M_I$ by Filon-Simpson-integration of $F_m$ in (\ref{eq.Fhw}) over
the interval $0\le y\le 1$ (Appendix \ref{app.FS}) and
\begin{equation}
M_I = \int_1^m e^{i\pi x}x^{1/x}dx +\frac{i}{\pi}\left[(-)^m m^{1/m}-1\right]
+\frac{i}{\pi}F_m;\quad g(x)=x^{1/x}\frac{1-\log x}{x^2};\quad (m=3).
\end{equation}
\end{algo}
With this choice, $g(x)$ has a maximum near $x=4.3$, associated
with the zero of $f^{(2)}(x)$, which (after detailed inspection) does
not destroy the alternating property---if higher order $f^{(n)}(x)$ were
employed, $m$ would have to be chosen differently.

\begin{table}
\caption{Convergence of Algorithm \ref{alg.as} with a Simpson integration on
$n$ abscissa points over $[0,1]$,
truncating the alternating series after the $l$-th term followed by extrapolation \cite{CohenExpM9}.}
\begin{tabular}{rrrr}
$n$ & $l$ & $\Re(M_I)$ & $\Im(M_I)$
\\
\hline
32&60&0.07077603721021819390&-0.68400038753980049654
\\
32&70&0.07077603721021819390&-0.68400038753980049654
\\
64&70&0.07077603918043104863&-0.68400038931936816305
\\
128&70&0.07077603930333884842&-0.68400038943052296458
\\
256&70&0.07077603931101698843&-0.68400038943746907333
\\
512&70&0.07077603931149681599&-0.68400038943790318846
\\
1024&70&0.07077603931152680433&-0.68400038943793032039
\\
2048&70&0.07077603931152867859&-0.68400038943793201613
\\
4096&70&0.07077603931152879573&-0.68400038943793212212
\\
8192&50&0.07077603931152880305&-0.68400038943793212874
\\
8192&60&0.07077603931152880305&-0.68400038943793212874
\\
16384&50&0.07077603931152880351&-0.68400038943793212916\\ 
16384&60&0.07077603931152880351&-0.68400038943793212916\\ 
\hline
\end{tabular}
\label{tab.alg2}
\end{table}
The speed of convergence of the method is demonstrated with Table \ref{tab.alg2}.
As the number of evaluations of $g$ is the product of
$n$ and $l$, it turns out effectively to be slower than Algorithm \ref{alg.esc}.

\subsection{Taylor Series the Logarithmic Term} \label{sec.Tlog}
The most advanced method expands the non-oscillatory term of (\ref{eq.ingr}) into
the Taylor series of the exponential:
\begin{algo} \label{alg.FH}
(Logarithmic expansion of $x^{1/x}$)
\begin{multline}
\int_1^{2N} e^{i\pi x}x^{1/x}dx = \int_1^{2N} e^{i\pi x}e^{\frac{1}{x}\log x}dx
= \int_1^{2N} e^{i\pi x}[1+\sum_{n\ge 1}\frac{1}{n!}\frac{\log ^nx}{x^n}] dx
\\
= -\frac{2i}{\pi}+\sum_{n\ge 1}\frac{1}{n!} \int_1^{2N} e^{i\pi x}\frac{\log ^nx}{x^n} dx
.
\end{multline}
\end{algo}
The integrals are the $V(\pi,n,n)$ defined in
equation (\ref{eq.Vksdef}) in Appendix \ref{app.FH}.
Summing over values taken from
Tables \ref{tab.ImV1s}--\ref{tab.Vks}
produces Table \ref{tab.FH}. 
\begin{table}
\caption{
Convergence of the partial sum of Algorithm \ref{alg.FH}.
}
\begin{tabular}{rrr}
$\max n$ & $\Re(M_I)$ & $\Im(M_I)$ \\
\hline
1 & 0.05762490298863188764 & -0.68331060191932132015 \\
2 & 0.06935454902524610824 & -0.68451362283943263006 \\
3 & 0.07066781734932533318 & -0.68408393964817446557 \\
4 & 0.07076978736326279015 & -0.68400839361204835470 \\
5 & 0.07077575770475264223 & -0.68400096184084805332 \\
6 & 0.07077602957520227166 & -0.68400042266805241181 \\
7 & 0.07077603908225272182 & -0.68400039107408096452\\
8 & 0.07077603931050590528 & -0.68400038950815469850\\
9 & 0.07077603931177817100 & -0.68400038944060977421\\
\hline
\end{tabular}
\label{tab.FH}
\end{table}
An accuracy of $10^{-19}$ in real and imaginary part of $M_I$ requires
summation up to $n=15$---an estimation derived from results
of Remark \ref{rem.tele}---, and has not been worked out.

\subsection{Contour Deformation} \label{sec.cont}
The path of the integration
may be deformed to a straight line (hypotenuse) from $z=1$ to $z=2N(1+\tau i)$
with adjustable slope $\tau>0$ towards the real axis,
and a straight line (short leg) parallel to the imaginary axis back to $2N$
on the real axis.
The contribution of $\Im z$ to $\exp(i\pi z)$ leads to 
a exponential reduction of the integrand as the distance to the real
axis grows. (More complicated paths appear not to be more efficient.)
The short leg of this triangle contributes
\begin{equation}
\lim_{N\to \infty}\int_{2N(1+\tau i)}^{2N} e^{i\pi z +\log z/z} dz
= -\frac{i}{\pi}
\end{equation}
to the integral.
\begin{algo} \label{alg.cont}
Compute 
\begin{equation}
M_I = -\frac{i}{\pi}+ \lim_{N\to \infty}\int_1^{2N(1+\tau i)} e^{i\pi z+\log z/z}dz
\end{equation}
with a Simpson integration on $s$ points $z_j=1+j(1+\tau i)\Delta t$,
$j=0,\ldots s$.
\end{algo}
\begin{table}
\caption{Algorithm \ref{alg.cont}: Simpson integration along the line from $z=1$ to
$z=2N(1+\tau i)$ with $s$ abscissa points.
}
\begin{tabular}{rrlrr}
$s$ & $N$ & $\tau$ & $\Re(M_I)$ & $\Im(M_I)$ \\
\hline
64000 & 20 & 0.1&0.07077602653194263984 & -0.68400205509417633076\\
64000 & 20 & 0.2 & 0.07077603931143220832 & -0.68400038944591701234\\
64000 & 20 & 0.3 & 0.07077603931156096762 & -0.68400038943797230450\\
128000 & 20 & 0.3 & 0.07077603931153081310 & -0.68400038943793467559\\
256000 & 20 & 0.3 & 0.07077603931152892845 & -0.68400038943793232378\\
\hline
64000 & 40 & 0.1&0.07077603931142945824 & -0.68400038944426949930\\
128000 & 40 & 0.1&0.07077603931149682353 & -0.68400038944357500783\\
64000 & 40 & 0.2 & 0.07077603931176127606 & -0.68400038943868543692\\
128000 & 40 & 0.2 & 0.07077603931154333296 & -0.68400038943797921085\\
\hline
64000 & 80 & 0.1&0.07077603931041311002 & -0.68400038945008174688\\
128000 & 80 & 0.1&0.07077603931145906913 & -0.68400038943869147372\\
256000 & 80 & 0.1&0.07077603931152444508 & -0.68400038943797958811\\
64000 & 80 & 0.2 & 0.07077603931534354347 & -0.68400038945029316088\\
128000 & 80 & 0.2 & 0.07077603931176721758 & -0.68400038943870468984\\
256000 & 80 & 0.2 & 0.07077603931154370430 & -0.68400038943798041416\\
512000 & 80 & 0.2 & 0.07077603931152973483 & -0.68400038943793514699\\
1024000 & 80 & 0.2 & 0.07077603931152886174 & -0.68400038943793231779\\
\hline
128000 & 160 & 0.1&0.07077603931039901029 & -0.68400038945023529554\\
256000 & 160 & 0.1&0.07077603931145818783 & -0.68400038943870107039\\
128000 & 160 & 0.15 & 0.07077603931285988807 & -0.68400038945052460466\\
256000 & 160 & 0.15 & 0.07077603931161199070 & -0.68400038943871915335\\
512000 & 160 & 0.15 & 0.07077603931153400264 & -0.68400038943798131810\\
\hline
\end{tabular}
\label{tab.cont}
\end{table}
The dependence on three configuration parameters makes quality assessment
more difficult than with the other methods, illustrated by Table \ref{tab.cont}\@.
The contribution missing for any finite $N$ is estimated by
\begin{equation}
\int_{2N(1+\tau i)}^\infty e^{i\pi z+\log z/z}dz
\approx
\int_{2N(1+\tau i)}^\infty e^{i\pi z}dz
=
\frac{i}{\pi}e^{2\pi i N(1+i\tau)},
\end{equation}
which serves as a guideline how large $\tau$ ought be made given a targeted
accuracy and an upper limit $N$.

\section{Summary} 
The value of $M_I$ is
\begin{multline}
\lim_{N\to \infty} \int_1^{2N} e^{i\pi x}x^{1/x} dx
\\
\approx
0.0707760393115288035395-0.68400038943793212918i
.
\end{multline}
The integral features highly oscillatory behavior and a logarithmic factor
in the main integrand, which sets up an interesting test case
outside the range of methods that assume ``nice'' analytic properties
in the complex plane.

\appendix

\section{Filon-Simpson Rule}\label{app.FS}
The Simpson rule of integration of a function $G(y)$ is an interpolation
between three abscissa points
$(y_0,G(0))$, $(y_1,G(1))$ and $(y_2, G(2))$ by a quadratic polynomial
in $y$ \cite[(25.2.1)]{AS},
\[
G_m(y)=
\frac{(y-y_1)(y-y_2)}{(y_0-y_1)(y_0-y_2)}G(0)
+
\frac{(y-y_0)(y-y_2)}{(y_1-y_0)(y_1-y_2)}G(1)
+
\frac{(y-y_0)(y-y_1)}{(y_2-y_0)(y_2-y_1)}G(2)
\]
followed by integration of this polynomial in the limits $y_0\le y\le y_2$.
A refinement in our case, based on the same quadratic interpolation, is
\begin{multline}
\int_{y_0}^{y_0+2h} e^{i\pi y}G_m(y)dy=
\frac{i e^{y_0\pi i}}{2h^2\pi^3}[-2+2h^2\pi^2-3ih\pi+(2-h\pi i)e^{2i\pi h}]G(0)
\\
-\frac{2 e^{y_0\pi i}}{h^2\pi^3}[-i+h\pi+(i+h\pi)e^{2i\pi h}]G(1)
\\
-\frac{i e^{y_0\pi i}}{2h^2\pi^3}[2+h\pi i+(-2+2h^2\pi^2+3i h\pi)e^{2i\pi h}]G(2),
\end{multline}
supposing equidistant abscissa points $y_1=y_0+h$ and $y_2=y_0+2h$.
This explicit recognition of the exponential factor
gains roughly one additional digit in accuracy in Table \ref{tab.alg2}
compared to the evaluation with $e^{i\pi y}$ incorporated in the
value of $G$.

\section{Fichtenholz Integrals}\label{app.FH}
\subsection{Fundamental Form}\label{app.V11}
In this appendix, a set of integrals $V(a,k,s)$ is targeted as an aid
to Algorithm \ref{alg.FH}.
This is basically evaluating
the generalized integro-exponential function \cite{MilgramMatC44,KolbigMathComp41} for a complex-valued parameter $z=-ia$.
\begin{defn} (Generalized Integro-Exponential Function)
\begin{equation} 
V(a,k,s)\equiv \int_1^\infty e^{iax}\frac{\log^k x}{x^s}dx
.
\label{eq.Vksdef}
\end{equation}
\end{defn}
\begin{rem}
There are three simple extensions:
\begin{itemize}
\item
Cases with a scaling factor $b$ in the logarithm can be reduced
to this fundamental form by binomial expansion of $(\log b+\log x)^k$,
\begin{equation}
\int_1^\infty e^{iax}\frac{\log^k(bx)}{x^s}dx
=
\sum_{l=0}^k \binom{k}{l} \log^{k-l}(b)\, V(a,l,s)
.
\end{equation}
\item
Reading
\begin{equation}
\int_1^\infty e^{iax}\frac{\log^k(bx)}{x^s}dx
=b^{s-1}\int_{b}^\infty e^{iax/b} \frac{\log^k y}{y^s}dy
,
\end{equation}
(obtained through the substitution $bx\to y$) from right to left
shows that other lower limits than 1 are also accessible once the $V$ are known
for general $a$.
\item
Integer powers of sines or cosines at the place of the exponential
lead back to the fundamental form via
Euler's formula:
\begin{eqnarray}
\int_1^\infty \sin^m(ax)\frac{\log ^k x}{x^s}dx
=
\frac{1}{(2i)^m}\sum_{l=0}^m\binom{m}{l}(-1)^l V[a(m-2l),k,s];
\\
\int_1^\infty \cos^m(ax)\frac{\log ^k x}{x^s}dx
=
\frac{1}{2^m}\sum_{l=0}^m\binom{m}{l} V[a(m-2l),k,s].
\end{eqnarray}
\end{itemize}
\end{rem}

Real and imaginary part of the value $V(a,1,1)$
\begin{equation}
\int_1^\infty e^{ia x}\frac{\log x}{x}dx =
\int_1^\infty \cos(a x)\frac{\log x}{x}dx
+i
\int_1^\infty \sin(a x)\frac{\log x}{x}dx
\end{equation}
are computed separately. The imaginary part is
\begin{equation}
\int_1^\infty \sin(a x)\frac{\log x}{x}dx
=
\int_0^\infty \sin(a x)\frac{\log x}{x}dx
-
\int_0^1 \sin(a x)\frac{\log x}{x}dx
\label{eq.sinFH}
\end{equation}
with one constituent \cite[(4.421.1)]{GR}\cite[(865.63)]{Dwight}\cite{Fichtenholz}
\begin{equation}
\int_0^\infty \sin(a x)\frac{\log x}{x}dx = -\frac{\pi}{2}(\gamma+\log a)
\approx -2.704825746060380848849568\quad (a=\pi)
\label{eq.sin0inf}
.
\end{equation}
The integral over $[0,1]$  is evaluated by Taylor expansion of the sine
which leads to the well converging representation
\begin{multline}
\int_0^1 \sin(a x)\frac{\log x}{x}dx
= -a\sum_{n\ge 0}(-1)^n \frac{a^{2n}}{(2n+1)!(1+2n)^2}
\\
= -a\,_2F_3\left(\begin{array}{c}1/2,1/2\\ 3/2,3/2,3/2
\end{array}
\mid-\frac{a^2}{4}
\right)
\approx -2.6581349165086 \quad (a=\pi)
.
\end{multline}
The difference between this value and (\ref{eq.sin0inf})
represents (\ref{eq.sinFH}),
\begin{multline}
\int_1^\infty \sin(ax)\frac{\log x}{x}dx
=
-\frac{\pi}{2}(\gamma+\log a)+a\sum_{n\ge 0}(-1)^n\frac{a^{2n}}{(2n+1)!(1+2n)^2}
\label{eq.ImV11}
\end{multline}
The value at $a=\pi$ is the head entry in Table \ref{tab.ImV1s}.

\begin{table}
\caption{Table of the imaginary part of $V(\pi,1,s)$.
}
\begin{tabular}{rr}
$s$ & $\Im V(\pi,1,s)$ \\
\hline
1 & -0.046690829551739977074516092264 \\
2 & -0.050400599397438879223041567776 \\
3 & -0.044723677797644192936988882199 \\
4 & -0.036181141258573216997609321919 \\
5 & -0.027931519676734642467423612590 \\
6 & -0.021096986691682143229642825070 \\
\hline
\end{tabular}
\label{tab.ImV1s}
\end{table}

The real part is started from \cite[(3.761.9)]{GR}
\begin{equation}
\int_0^\infty \frac{\cos (a x)}{x^{1-\mu}}dx = \frac{\Gamma(\mu)}{a^\mu}
\cos(\frac{\mu\pi}{2})
\label{eq.V01}
,
\end{equation}
which is differentiated with respect to $\mu$ with the product rule,
\begin{equation}
\int_0^\infty \cos (a x)\frac{\log x}{x^{1-\mu}} dx=
\frac{\Gamma(\mu)}{a^\mu}
\left(
\psi(\mu) \cos\frac{\mu\pi}{2}
-\ln a \cos\frac{\mu\pi}{2}
-\frac{\pi}{2} \sin\frac{\mu\pi}{2}
\right)
.
\label{eq.cosmu}
\end{equation}
The reflection formulas for the $\Gamma$-function and Digamma function $\psi$
show that in the limit $\mu\to 0$ \cite[(6.1.17),(6.3.7)]{AS}
\begin{multline}
\int_0^\infty \cos (a x)\frac{\log x}{x^{1-\mu}} dx \to
-\frac{1}{\mu^2}-\frac{\pi^2}{24}
+\gamma\left(\frac{\gamma}{2}+\ln a\right)+\frac{\ln^2 a}{2}
+O(\mu)
\label{eq.cos0}
\end{multline}
where $\gamma\approx 0.57721$ is the Euler-Mascheroni
constant \cite[A001620]{EIS}. Expansion of the cosine
in the standard Taylor series \cite[(1.411.3)]{GR} and interchange of integration and summation yields
\begin{equation}
\int_0^1 \frac{\cos(a x)}{x^{1-\mu}}dx
=\sum_{n\ge 0}\frac{(-a^2)^n}{(2n)!(2n+\mu)}.
\end{equation}
Differentiation with respect to $\mu$ introduces the logarithm,
\begin{equation}
\int_0^1 \cos(a x)\frac{\log x}{x^{1-\mu}}dx
= -\sum_{n\ge 0}\frac{(-a^2)^n}{(2n)!(2n+\mu)^2}.
\end{equation}
We subtract this from (\ref{eq.cos0}) and notice that the $\sim 1/\mu^2$-singularity
cancels with the term $n=0$ of the series as $\mu\to 0$,
\begin{equation}
\int_1^\infty
\cos(a x)\frac{\log x}{x}dx
=
-\frac{\pi^2}{24}
+\gamma\left(\frac{\gamma}{2}+\ln a\right)+\frac{\ln^2 a}{2}
+\sum_{n\ge 1}\frac{(-a^2)^n}{(2n)!(2n)^2}.
\end{equation}
The sum
converges quickly,
the value at $a=\pi$ is the first entry in Table \ref{tab.ReV1s}.

\begin{table}
\caption{Table of the real part of $V(\pi,1,s)$.
}
\begin{tabular}{rr}
$s$ & $\Re V(\pi,1,s)$ \\
\hline
1 & 0.057624902988631887643485542240 \\
2 & 0.029913203983934978439301792236 \\
3 & 0.010937363639874260291206201403 \\
4 & -0.000250069139610211209861368961 \\
5 & -0.006024230915536561502482189260 \\
6 & -0.008508918812024751462009533761 \\
\hline
\end{tabular}
\label{tab.ReV1s}
\end{table}

Combined with (\ref{eq.ImV11}) this reads
\begin{equation}
\int_1^\infty e^{iax}\frac{\log x}{x}dx
=
-\frac{\pi^2}{24}
+\gamma\left(\frac{\gamma}{2}+\ln a\right)+\frac{\ln^2 a}{2}
-\frac{\pi}{2}i(\gamma+\log a)
+\sum_{n\ge 1}\frac{(ia)^n}{n!n^2},
\label{eq.coFHa}
\end{equation}
which constitutes a ``root'' value $V(a,1,1)$ in the tree of integrals (\ref{eq.Vksdef}).
A summary of (\ref{eq.cosmu}) and the associated imaginary part is
\begin{equation}
\int_1^\infty e^{iax}\frac{\log x}{x^{1-\mu}}dx = V(a,1,1-\mu)=
\frac{\Gamma(\mu)}{a^\mu}e^{i\pi\mu/2}\left(\psi(\mu)-\ln a +\frac{i\pi}{2}\right)
+\sum_{n\ge 0}\frac{(ia)^n}{n!(n+\mu)^2}
.
\end{equation}
\subsection{Higher Powers of  the Rational}
Integration of (\ref{eq.coFHa}) with respect to $a$
is a measure to increase the parameter $s$, the power of $x$ in the denominator,
by one:
\begin{equation}
\int_0^a da'
\int_1^\infty dx
e^{ia'x}\frac{\log x}{x}
=
\int_1^\infty
[e^{iax}-1]\frac{\log x}{ix^2}dx
=
\frac{1}{i}\int_1^\infty
e^{iax}\frac{\log x}{x^2}dx
-\frac{1}{i}
.
\label{eq.Vinta}
\end{equation}
\begin{multline}
\int_1^\infty
e^{iax}\frac{\log x}{x^2}dx
=
1+ia
\left[
-\frac{\pi^2}{24}+\gamma(\frac{\gamma}{2}+\ln a -1)
+\frac{\ln^2a}{2}-\ln a +1 -\frac{i \pi}{2}(\gamma+\ln a-1)
\right]
\\
+\sum_{n\ge 1}\frac{(ia)^{n+1}}{(n+1)!n^2}
,
\end{multline}
where \cite{RoyAMM94}
\begin{equation}
\sum_{n\ge 1}\frac{(ia)^{n+j}}{(n+j)!n^2} = \frac{(ai)^{1+j}}{(1+j)!}
\,_3F_3\left(\begin{array}{c}1,1,1\\2,2,2+j\end{array}\mid ai\right)
.
\end{equation}
Inserting $a=\pi$
yields the second lines in Table \ref{tab.ReV1s} and \ref{tab.ImV1s}.
The concept of (\ref{eq.Vinta}) generalizes to higher powers $k$ of the logarithm,
\begin{equation}
\int_0^a da' V(a',k,s) = \frac{1}{i}\int_1^\infty [e^{iax}-1]\frac{\log^k x}{x^{1+s}}dx
=
-iV(a,k,1+s) +iV(0,k,1+s)
.
\label{eq.Vintagen}
\end{equation}
The last term, the non-oscillatory integrals at $a=0$, are known \cite[(4.272.6)]{GR},
\begin{equation}
V(0,0,s) = \frac{1}{s-1},\quad
V(0,k,s) = \frac{\Gamma(k+1)}{(s-1)^{k+1}},
\end{equation}
Milgram's equation (2.29) \cite{MilgramMatC44}.
Iterated application of this rule computes the
chain of $V(a,1,1)\to V(a,1,2) \to V(a,1,3)\to\ldots$ as follows:
\begin{multline}
\int_1^\infty
e^{iax}\frac{\log x}{x^3}dx
=
\frac{\pi^2a^2}{48}-\frac{\gamma a}{4}
(\gamma a+2a \ln a-3a)
-\frac{a^2}{4}(\ln^2 a -3\ln a+7/2)
+\frac{1}{4}
\\
+\frac{\pi a^2}{4}i(\gamma+\ln-\frac{3}{2})
+ai
+\sum_{n\ge 1}\frac{(ia)^{n+2}}{(n+2)!n^2}
.
\end{multline}

\begin{multline}
\int_1^\infty
e^{iax}\frac{\log x}{x^4}dx
=
\frac{\pi^2 a^3}{144}i
-\frac{\gamma a^3}{12}i
(\gamma+2\ln a-11/3)
+\frac{a^3}{3}i(-\frac{85}{72}+\frac{11}{12}\ln a-\frac{1}{4}\ln^2 a)
+\frac{a}{4}i
\\
-\frac{\pi a^3}{12}(\gamma+\ln a-11/6)-\frac{a^2}{2}+\frac{1}{9}
+\sum_{n\ge 1}\frac{(ia)^{n+3}}{(n+3)!n^2}
.
\end{multline}
The cases of $a=\pi$ are in Tables \ref{tab.ImV1s} and \ref{tab.ReV1s}.

\subsection{Higher Powers of the Logarithm}\label{app.logp}

As seen in Section \ref{app.V11}, differentiation with respect to the parameter $\mu$
increases the power of the logarithm:
\begin{equation}
\frac{d}{d s}V(a,k,s) = -\int_1^\infty e^{iax}\frac{\ln^{k+1}x}{x^s}dx
= -V(a,k+1,s)
.
\label{eq.dVds}
\end{equation}

To support differentiation with respect to the $s$-parameter, we recompute
$V(a,1,2-\mu)$, which we re-integrate over $a$ as in (\ref{eq.cosmu}):
\begin{multline}
\int_0^a da' \int_0^\infty \cos (a' x)\frac{\log x}{x^{1-\mu}} dx=
\int_0^\infty \sin (a x)\frac{\log x}{x^{2-\mu}} dx
\\
=
\frac{\Gamma(\mu)}{a^{\mu-1}(1-\mu)}
\left[ \left( \psi(\mu) +\frac{1}{1-\mu} -\ln a\right) \cos\frac{\mu\pi}{2} -\frac{\pi}{2} \sin\frac{\mu\pi}{2} \right]
.
\end{multline}
The integration is applied in parallel to the complementary interval $0\le x\le 1$,
\begin{equation}
\int_0^1 \sin (a x)\frac{\log x}{x^{2-\mu}} dx=
-\sum_{n\ge 0} \frac{(-)^na^{2n+1}}{(2n+1)!(2n+\mu)^2}
,
\end{equation}
and the difference is
\begin{multline}
\int_1^\infty \sin (a x)\frac{\log x}{x^{2-\mu}} dx
=
\Im V(a,1,2-\mu)
\\
=
\frac{\Gamma(\mu)}{a^{\mu-1}(1-\mu)}
\left[ \left( \psi(\mu) +\frac{1}{1-\mu} -\ln a\right) \cos\frac{\mu\pi}{2} -\frac{\pi}{2} \sin\frac{\mu\pi}{2} \right]
+\sum_{n\ge 0} \frac{(-)^na^{2n+1}}{(2n+1)!(2n+\mu)^2}
.
\end{multline}
This is differentiated with respect to $\mu$, and the singularities
$\sim 2a/\mu^3$ from the incomplete Gamma-function and the $n$-sum cancel
in the limit $\mu\to 0$:
\begin{multline}
\int_1^\infty \sin (a x)\frac{\log^2 x}{x^2} dx= a\Big[
\gamma(-\gamma+2-\ln a)
\ln a
+\frac{1}{12} \pi^2(\ln a-1+\gamma)
-\frac{2}{3}\zeta(3)
\\
-2\ln a
+\ln^2 a-\frac{1}{3}\ln^3 a
-\frac{\gamma^3}{3}
-2\gamma+\gamma^2
+2
\Big]
-2a\sum_{n\ge 1} \frac{(-a^2)^n}{(2n+1)!(2n)^3}
.
\end{multline}
The explicit value at $a=\pi$ is
\begin{equation}
\int_1^\infty \sin (\pi x)\frac{\log^2 x}{x^2} dx=
\Im V(\pi,2,2) \approx
-0.00240604184022261982751961704408
.
\end{equation}

This demonstrates the technique.
Starting from $V(a,2,2)$, a ladder of integrals is constructed 
according to (\ref{eq.Vintagen}). Each time, a term $\sim (-1)^k k!(ia)^{s-1}/\mu^{k+1}$
cancels---carried over from the simple pole of the $\Gamma$-function (\ref{eq.V01})
through $k$ differentiations and $s-1$ integrations---when
the complementary integrals of $0\le x< \infty$ and $0\le x\le 1$
are combined.
Numerical examples are gathered in Table \ref{tab.Vks}.
\begin{table}
\caption{Table of $V(\pi,k,s)$.
}
\begin{tabular}{rrrr}
$k$ & $s$ & $\Re V(\pi,k,s)$ & $\Im V(\pi,k,s)$ \\
\hline
2 & 2 & 0.0234592920732284411947929 & -0.0024060418402226198275196\\
2 & 3 & 0.0147167653246107850628908 & -0.0078739079798083225909317\\
2 & 4 & 0.0080788243950624846223234 & -0.0087094002295813942921939\\
2 & 5 & 0.0038282314868382609838401 & -0.0076206659960747444416239\\
2 & 6 & 0.0013846381369360916659413 & -0.0060334458936144149013900\\
2 & 7 & 0.0000998710536555974796592 & -0.0045501919692512628703906\\
2 & 8 & -0.0005097889780015996357674 & -0.0033548439064072548162349\\
\hline
3 & 3 & 0.0078796099444753496396155 & 0.0025780991475489869676813\\
3 & 4 & 0.0053790452810071493354663 & -0.0004578919913424486140616\\
3 & 5 & 0.0032308011441634551812315 & -0.0014908072125213319989936\\
3 & 6 & 0.0017674846795167938450865 & -0.0015900953082259670812147\\
3 & 7 & 0.0008825074832995147577407 & -0.0013496431705418864940445\\
3 & 8 & 0.0003872374479338897353530 & -0.0010417218132785702486827\\
\hline
4 & 4 & 0.0024472803344989672143380 & 0.0018131048670266608466124\\
4 & 5 & 0.0018067919115537410157224 & 0.0004312822675128498086469\\
4 & 6 & 0.0011430051029650516946321 & -0.0001368354074002352164914\\
4 & 7 & 0.0006599851739729340692718 & -0.0003012860412794823958833\\
4 & 8 & 0.0003564954008069045751359 & -0.0002990689541546756726797\\
\hline
5 & 5 & 0.0007164409787822497618537 & 0.0008918125440361658853376\\
5 & 6 & 0.0005826627556204031506092 & 0.0003133177357343641402474\\
5 & 7 & 0.0003859348621737006202841 & 0.0000540098043633385948340\\
5 & 8 & 0.0002304000284892623408819 & -0.0000404135204291844446708\\
\hline
6 & 6 & 0.0001957467237331888336746 & 0.0003882044128618852155461 \\
6 & 7 & 0.0001826715069173610144948 & 0.0001565025492374147392079 \\
6 & 8 & 0.0001272475588547481432200 & 0.0000473426202252359704144 \\
\hline
\end{tabular}
\label{tab.Vks}
\end{table}

\begin{rem} \label{rem.tele}
Partial integration
\begin{multline}
\int e^{iax}\frac{\log ^k x}{x^s}dx
=
\int e^{iax}\frac{\log ^k x}{x^{s-1}}\frac{1}{x}dx
\\
=
e^{iax}\frac{\log ^{k+1} x}{x^{s-1}}
-
\int \left[
ia e^{iax}\frac{\log ^k x}{x^{s-1}}
+ e^{iax}k \frac{\log ^{k-1} x}{x^s}
+ e^{iax} (1-s)\frac{\log ^k x}{x^s}
\right]\log x dx
\end{multline}
 yields a contiguous relation for three points in a triangle
in the square grid of $(k,s)$-pairs (Milgram's equation (2.4) \cite{MilgramMatC44}):
\begin{equation}
V(a,k,s)
=
- \frac{ia}{1+k} V(a,k+1,s-1)
+ \frac{s-1}{1+k}V(a,k+1,s)
.
\label{eq.contig}
\end{equation}
Since the calculation of values at large $k$ is a laborious task,
this formula offers a route to cheaper
calculation by (i) tabulation of $V(a,k,s)$ at some small $k$ up to a rather
large $\hat s$, involving only integrations of powers of $\ln a$ multiplied
by powers of $a$ \cite[(2.722)]{GR},
(ii) numerical calculation of the $V(a,k+\Delta k,\hat s)$
up to the desired $k$ with some brute force method like the $1/x$ mapping (\ref{eq.1u}), which 
converges well since $\hat s$ is large, (iii) telescoping from $\hat s$
backwards with (\ref{eq.contig}) to fill the table for increasing $k$ and
decreasing $s$.
\end{rem}
\begin{rem}
The last term in (\ref{eq.contig}) can be eliminated via (\ref{eq.dVds}),
\begin{equation}
(1+k)V(a,k,s)
=
- ia V(a,k+1,s-1)
+ (1-s)\frac{d}{ds}V(a,k,s)
,
\end{equation}
The solution of this inhomogeneous differential equation
is
\begin{equation}
V(a,k,s) = \frac{1}{(s-1)^{1+k}}\left(-ia\int V(a,k+1,s-1)(s-1)^kds +const \right).
\end{equation}
\end{rem}

\section{Inverse Euler-Maclaurin} \label{app.EMcl} 
A standard idea of integration is to split the integral over intervals
that are commensurable
with the frequency of the oscillation,
to replace the generic factor $g$ in the integral by some approximation which
allows integration in closed form---assuming that a massive cancellation
can be obtained---, then to gather the sum over the intervals
with some Euler-Maclaurin approach \cite{BerndtAA28}.
Applied to (\ref{eq.Fhw}),
\begin{equation}
F_m\equiv \int_m^\infty e^{i\pi x}g(x) dx
=\sum_{k=m+1,m+3,m+5,\ldots}^\infty \int_{k-1}^{k+1} e^{i\pi x}g(x)dx,
\end{equation}
$g$ is approximated
by its Taylor series,
\cite{MatharArxiv0902,GouldAMM70},
\begin{equation}
g(x)=\sum_{d=0}^\infty \frac{1}{d!}g^{(d)}(k)(x-k)^d
.
\end{equation}
\begin{equation}
F_m=\sum_{k=m+1,m+3,m+5,\ldots}^\infty e^{i\pi k}\int_{-1}^{1}dy e^{i\pi y}\sum_{d=0}^\infty \frac{1}{d!}g^{(d)}(k)y^d
.
\label{eq.Fmm}
\end{equation}
\begin{defn} (Moments of Filon Quadratures)
\begin{equation}
S(d)
\equiv \int_{-1}^{1} e^{i\pi y}y^d dy; \quad \bar S(d)\equiv \frac{S(d)}{d!};\quad d=0,1,2,3,\ldots
\end{equation}
\end{defn}
Initial values are
\begin{equation}
S(0)=0;\quad S(1)=\frac{2i}{\pi}.
\end{equation}
The recurrence is
\begin{equation}
S(d+1)=\frac{i}{\pi}\left\{
\left[1+(-1)^d\right]+(d+1)S(d)
\right\}
.
\end{equation}
With \cite[(3.761.5),(3.761.10)]{GR} the cases for even and odd indices are
\begin{eqnarray}
\bar S(2d+1) &=& 2i\frac{(-)^d}{\pi^{2d+1}}\sum_{j=0}^d \frac{1}{(2j+1)!}(-\pi^2)^j
.\\
\bar S(2d) &=& 2\frac{(-)^d}{\pi^{2d}}\sum_{j=0}^{d-1} \frac{1}{(2j+1)!}(-\pi^2)^j
.
\end{eqnarray}
This rewrites (\ref{eq.Fmm}),
\begin{equation}
F_m=
\sum_{k=m+1,m+3,m+5,\ldots}^\infty e^{i\pi k} \sum_{d=0}^\infty \bar S(d)g^{(d)}(k)
=
(-1)^{m+1}
\sum_{l=0}^\infty \sum_{d=1}^\infty \bar S(d)g^{(d)}(m+2l+1)
.
\end{equation}
The Euler-Maclaurin formula proposes to replace the sum over the $d$-th derivatives
by
\begin{multline}
\sum_{l=0}^\infty g^{(d)}(m+2l+1)
=
\frac{1}{2}\bigg(
g^{(d)}(m+1)+\int_{m+1}^\infty g^{(d)}(x)dx
\\
-\sum_{D=1,3,5\ldots}
\frac{B_{D+1}}{(D+1)!}2^{D+1}g^{(d+D)}(m+1)
\bigg)
.
\end{multline}
$F_m$ becomes a double sum over $d$ and $D$, and resummation
proposes Algorithm \ref{alg.em} to accumulate the derivatives of
the (smooth) function $g$ to calculate
\begin{equation}
M_I=\frac{i}{\pi}(F-2); \quad g\equiv x^{1/x}\frac{1-\log x}{x^2}.
\end{equation}
\begin{algo} \label{alg.em}
(1-sided Fourier-Euler-Maclaurin)
\begin{multline}
2(-)^{m+1}F
=
-\bar S(1)g^{(0)}(m+1) \\
+
\sum_{d=1}^\infty
\left(
\bar S(d)
-
\bar S(d+1)
-
\sum_{l=0}^{\lfloor (d-1)/2\rfloor}
\bar S(d-1-2l)
\frac{B_{2+2l}}{(2+2l)!}2^{2+2l}
\right)
g^{(d)}(m+1)
.
\nonumber
\end{multline}
\end{algo}
Implementation of this approach reveals that the sum over the $d$-th derivatives
shows converging behavior only up to $d\approx 6$. I attribute this
to the same
logarithmic branch cut that
constraints the useful depths of the partial integrations in Algorithms
\ref{alg.pI} and \ref{alg.esc}.

\bibliographystyle{amsplain}
\bibliography{all}

\end{document}